\documentclass{amsart}

\title[Non-autonomous dynamics]{Non-autonomous Dynamics in $\PP^k$}
\author{Han Peters}
\date{November 14, 2003}

\usepackage{amscd}

\newtheorem{theorem}{Theorem}
\newtheorem{lemma}[theorem]{Lemma}
\newtheorem{proposition}[theorem]{Proposition}
\newtheorem{corollary}[theorem]{Corollary}

\theoremstyle{definition}
\newtheorem{defin}[theorem]{Definition}

\theoremstyle{remark}
\newtheorem{remark}[theorem]{Remark}

\newcommand{\NN}{\mathbb{N}}
\newcommand{\RR}{\mathbb{R}}
\newcommand{\CC}{\mathbb{C}}
\newcommand{\PP}{\mathbb{P}}

\begin{document}

\begin{abstract}
We study the dynamics of compositions of a sequence of holomorphic
mappings in $\PP^k$. We define ergodicity and mixing for
non-autonomous dynamical systems, and we construct totally
invariant measures for which our sequence satisfies these
properties.
\end{abstract}

\maketitle

\section{Introduction}

Non-autonomous dynamics differs from standard dynamics in that
instead of iterating a single map, we consider compositions of a
sequence of maps. The main goal of non-autonomous dynamics is to
generalize theorems that hold in the autonomous setting or to find
counterexamples. Here, we will try to generalize theorems which
state that for every
 complex mapping there exists a natural measure which is mixing
and thus ergodic. This was first proved by Brolin for polynomials
in the complex plane in [Br], by Bedford and Smillie for H\'{e}non
mappings in $\CC^2$ in [BS], and for regular polynomial mappings of
$\CC^k$ by Forn\ae ss and Sibony in [FS1]. It has also been shown for
endomorphisms of $\PP^k$, see for instance the articles by
Briend and Duval \cite{bd} or Guedj and Sibony \cite{gs}.

Non-autonomous systems of polynomials in the complex plane have
been studied in the past years by many authors, see for instance
the survey article by Comerford \cite{co}, which has an extensive
bibliography. It has turned out that a good setting in which to work
is that of bounded sequences of
monic polynomials of some fixed degree which were first
considered in \cite{fs2}. We will work with more general mappings, and our
results imply the same results for such sequences.

Non-autonomous systems in higher complex dimensions have been
studied only rarely. We will look at a compact sequence
 of holomorphic mappings on $\PP^k$, which we will define more
precisely in the next section. This setting has already been
studied in \cite{fw}, but in a rather different way. There, the
dynamics of all nearby mappings of a holomorphic mapping were
studied at the same time, while we will study the dynamics of one
fixed sequence.

Let $P_n$ be a compact sequence of holomorphic mappings, and let
$\mu_n$ be the equilibrium measures for this sequence, which we
will define later. The main results of this paper are the
following two theorems:

\begin{theorem}
The system $\left(\{P_n\},\{\mu_n \}\right)$ is randomly ergodic.
\end{theorem}

\begin{theorem}
The system $\left(\{P_n\} , \{\mu_n\}\right)$ is randomly mixing.
\end{theorem}

In Section 2 we will set our notation and give the precise
definitions of randomly ergodic and randomly mixing, and in the
Section 3 we will use pluripotential methods to  introduce the
equilibrium measures $\mu_n$. In Section 4 we will prove a
series of lemmas considering the convergence of preimages, and we
will give the proofs of our two main theorems in the
Section 5. In the last section we will prove that an autonomous
system which is randomly ergodic is in fact weakly mixing.

\section{Non-autonomous systems in $\PP^k$}

We will now introduce the setting for this paper. Let
$\mathcal{P}$ be a compact family ( in the coefficients topology)
of holomorphic endomorphisms of $\PP^k$ whose degrees are at least
$2$ and bounded from above, and let $P_0, P_1, \ldots$ be a sequence of polynomials in
$\mathcal{P}$, where $P_n$ has degree $d_n$.

We define
$$P(n) = P_n \circ \ldots P_1;\quad d(n)= d_n \cdot \ldots d_1$$
and for $n$ larger than $m$ we will write
$$
P(m,n) = P_n \circ \dots \circ P_{m+1};\quad  d(m,n) = d_n \cdot \ldots d_{m+1}
$$

For a point $z=z_0$ in $\PP^k$ we shall also write $z_n$ for
$P(n)(z)$, which we shall say is a point at stage $n$. Thus $P_n$
is a mapping from stage $n-1$ to stage $n$.

Recall that a measure preserving automorphism $f$ of a
space X with probability measure $\mu$ is called \emph{ergodic} if all
totally invariant measurable subsets $A$ of $X$ either have full
or empty measure, and that $f$ is called \emph{mixing} if for all measurable sets
$A$ and $B$ we have that
$$
\mu (f^{-n}(A)\cap B) - \mu(A)\cdot\mu(B) \rightarrow 0.
$$

We would like to study these two properties in the non-autonomous
setting, but in this setting the above definitions do not make much sense.
First of all, in general a sequence of maps $f_1, f_2, \ldots$
will not have a probability measure that is invariant for all
$f_n$. We shall say that $\{f_n\}$ is \emph{measure preserving} for a
sequence of probability measures $\mu_0, \mu_1, \ldots$, if
$f_{n*}\mu_{n-1} = \mu_n$ holds for every $n$. Secondly, there will
generally be no proper measurable subsets which are invariant for
all $f_n$. We say that a sequence $A_0, A_1, \ldots$ is \emph{totally
invariant} if $f_n^{-1}(A_n) = A_{n-1}$ for every $n$. We make
the following definitions.

\begin{defin}
A measure preserving sequence $\{f_n\}$ is \emph{randomly ergodic} if for
all totally invariant sequences $A_0, A_1, \ldots$, where $A_n$ is
a $\mu_n$-measurable set, we have that $\mu_n(A_n)$ is $0$ or $1$.
\end{defin}

\begin{defin}
A measure preserving sequence $\{f_n\}$ is \emph{randomly mixing} if for
all continuous functions $\phi$ and $\psi$ on $X$, we have that
$$
\int (\phi\circ f(n)) \cdot \psi \, d\mu_0 - \int \phi \, d\mu_n \int
\phi \, d\mu_0 \rightarrow 0.
$$
\end{defin}

Both definitions can also be studied in the autonomous setting,
where a single map is iterated. Since continuous functions are
dense in $\mathcal{L}^2(\mu)$, we have that randomly mixing and
mixing are equivalent. However, randomly ergodic is strictly
stronger property than ergodic. It is easy to check that
randomly ergodic implies ergodicity, but the only measure for
which an automorphism is randomly ergodic is a point mass at a
fixed point, which is certainly not the case for the classical
definition. We note that a randomly mixing system is not
necessarily randomly ergodic for the same reason.

It would be interesting to find a generalization of ergodicity that
is useful for the study of the dynamics of a sequence of automorphisms.

\section{Equilibrium Measures}

The following construction of the equilibrium measures is
fairly standard in holomorphic dynamics and can be found in \cite{si},
 and can also be found for non-autonomous systems in \cite{fw}.

Since $\mathcal{P}$ is compact, we can extend all mappings $P$
in $\mathcal{P}$ to homogeneous
polynomial mappings $\tilde{P}$ of $\CC^{k+1}$ in such a way that
the coefficients of every $\tilde P$ are bounded by some uniform
constant $M$, and such that the images of the unit sphere in
$\CC^{k+1}$ are bounded away from the origin. In other words, there
exists some constant $t>1$ such that
\begin{eqnarray} \label{convergence}
\frac{1}{t} \|z\|^{d_n} < \| \tilde{P_n}(z) \| < t \|z\|^{d_n},
\end{eqnarray}
holds for any nonzero $z$ in $\CC^{k+1}$ and any $n$.

For every $i \in \NN$ and $n\ge i$, we define the function
$$
G_{n,i}(z) := \frac{1}{d(i,i+n)} \log \| \tilde{P}(i,i+n)(z)\|.
$$

\begin{lemma}
As $n \rightarrow \infty$, the functions $G_{n,i}$ converge
uniformly on $\CC^{k+1}$ to a continuous and plurisubharmonic
function $G_i$.
\end{lemma}
\begin{proof}
Fix $\epsilon >0$. It follows from \eqref{convergence} that for
any $z$ in $\CC^{k+1}$ we have
\begin{eqnarray*}
|G_{n+1,i}(z) - G_{n,i}(z)| < \frac{\log(t)}{d(i, n+i+1)}.
\end{eqnarray*}
Therefore, we have for any $m \ge n$ that
\begin{eqnarray}
\label{green} |G_{m,i}(z) - G_{n,i}(z)| < \frac{\log(t)}{d(i, i+n)
(d_{n+1}-1)}.
\end{eqnarray}
Since every $d_n$ is at least $2$ we can choose $n$ large enough so that
$$|G_{m,i}(z) - G_{n,i}(z)| <\epsilon,$$
for any $m\ge n$. It follows that the sequence $G_{n,i}$ converges
uniformly to a limit map $G_i$, and since all the functions $G_{n,i}$ are
continuous and plurisubharmonic, the limit map is also
continuous and plurisubharmonic.
\end{proof}

It follows from \eqref{green} that $G(z) = \log\|z\| + O(1)$. Also,
since every $\tilde{P}_n$ is homogeneous, we have that
$G_i(\lambda z) = \log(\lambda) + G_i(z)$. We get the equation
\begin{eqnarray} \label{pullback}
\tilde{P}_n^* G_n = d_n G_{n-1}.
\end{eqnarray}

Let $\pi$ be the projection from $\CC^{k+1}$ to $\PP^k$. We can
define $(1,1)$ currents $T_i$ on $\PP^k$ which satisfy
\begin{eqnarray*}
\pi^* T_i := dd^c G_i.
\end{eqnarray*}
$T_i$ is a current of mass $1$ on $\PP^{k}$,
that does not depend on our choices for $\tilde{P}_n$. It follows
from equation
\eqref{pullback} that
$$P_n^*T_n = d_n T_n.$$

Since $G_n$ is continuous, it follows from \cite{bt} that
we can define $\mu_n = (T_n)^k$. Since $T_i$ has unit mass, we
get that $\mu_n$ is a probability measure and since $G_n$ is locally bounded
it follows from Proposition 4.6.4 in the book by Klimek \cite{kl}
that $\mu_n$ does not assign any mass to locally pluripolar sets.

We call $\mu_n$ the \emph{equilibrium measure} at stage $n$ and we have that
$P_n^*\mu_n = d_n^k \mu_{n-1}$, and that $P_{n*} \mu_{n-1} =
\mu_n$.

\section{Uniform Convergence of Preimages}

Recall the following theorem, which was proved by H. Brolin \cite{br}
for polynomials and by M. Lyubich \cite{ly} and independently by
A. Freire, A. Lopes and R. Ma\~n\'e \cite{flm} for rational functions:

\begin{theorem} \label{brolin}
Let $R(z)$ be a rational function of degree $d\ge 2$, and let $R^n$
be its n-th iterate. Then for all $a \in \PP\setminus \mathcal{E}_R,$
$\mathrm{card}(\mathcal{E}_R)\le 2$,
\begin{eqnarray*}
\frac{1}{d^n}(R^n)^*\delta_a \rightarrow \mu.
\end{eqnarray*}
\end{theorem}

Here $\delta_x$ is the dirac mass at $x$. It follows from Theorem
1.2 of \cite{rs} that this theorem can be generalized to our
setting. However, to prove Theorems 1 and 2 we will need the
uniform versions of this theorem which we will prove in this
section. Our proofs will be similar to the method used by Lyubich
to prove the above theorem, and which was later used by J. Briend
and J. Duval in \cite{bd} to prove similar results for
endomorphisms of $\PP^k$.

Define $\eta_{x,n,i}$ to be the probability measure with mass
$\frac{1}{d(i, n+i)^{k}}$ at all the preimages $P(i, n+i)^{-1}(x)$
counting multiplicity. In other words,

$$
\eta_{x,n, i} = \frac{P(i, n+i)^* \delta_x}{d(i,n+i)^{k}}.
$$
(For simplicity of notation, we shall write $\eta_{x, n}$ for
$\eta_{x, n, 0}$).

For two probability measures $\mu_1, \mu_2$ on $\PP^k$ we define the distance
$$
d(\mu_1, \mu_2) = \sup_\phi |\int \phi \, d\mu_1 - \int \phi \,
d\mu_2|,
$$
where the supremum is taken over all $\mathcal{C}^1(\PP^k)$
functions $\phi$ for which $|\phi(z)|$ and $|\nabla \phi(z)|$ are
bounded by $1$. It is clear that the topology induced by this distance
is weaker than the strong topology on probability measures. In fact
a sequence of probability
measures $\nu_n$ converges weakly to $\mu$ if and only if
$d(\nu_n, \mu) \rightarrow 0$ since we are working in a compact space.

The following proposition shows that as $n$ gets large, the measures
$\eta_{x,n}$ depend less and less on the point $x$.

\begin{proposition} \label{preimage}
Let $\epsilon > 0$. Then there exists an $N \in \NN$ and subsets $X_n$
of $\PP^k$ such that for every $n$ larger than $N$ we have that
$\mu_n(X_n) < \epsilon$, and also
$$
d(\eta_{n,x}, \eta_{n,y}) < \epsilon,
$$
for every $x,y$ outside of $X_n$.
\end{proposition}

The proof is given below.

Fix $\epsilon>0$, and let $l=l(\epsilon)$ be some large
enough number that we will define later. For $n$ greater
or equal to $l$, let $V_{l,n}$ be the set of critical
values of the holomorphic mapping $P(n-l,n)$.

\begin{lemma}
There exists a $\delta$ such that the $\mu_n$ mass of the
$\delta$-neighborhood of $V_{l,n}$ is less than $\epsilon$ for any
$n$ larger than $l$.
\end{lemma}

\begin{proof}
We have seen that the measures $\mu_n$ do not assign any mass
to pluripolar sets. Therefore, there exists for each $n \ge l$
a $\delta_n$ such that the $\delta_n$-neighborhood of $V_{l,n}$
has $\mu_n$ mass less than $\epsilon$. Let $\mathcal{S}$
be the set of sequences of polynomials of $\mathcal{P}$
with the product topology, so that $\mathcal{S}$ is a compact set.

The maps $G_{n,i}$ depend continuously on the sequence in
$\mathcal{S}$, and since $G_{n,i}$ converges uniformly to the map
$G_i$, we have that $G_i$ also depends continuously on $\mathcal{S}$.
Let $\{S^j\}$ be a sequence of sequences in $\mathcal{S}$ that
converges uniformly to $S \in \mathcal{S}$. Write $G_i^j, G_i$,
$\mu_i^j, \mu_i$ for the Green's functions and equilibrium measures
corresponding to the sequences $S^j$ and $S$. Then we have that
$G_i^j \rightarrow G_i$ uniformly on $\CC^k$, and therefore it follows
from \cite{cln} that $\mu_i^j$ converges weakly to $\mu_i$.

Since the sets of critical values $V_{l,n}$ also vary
continuously as a function in $\mathcal{S}$, we have
that $\delta_n$ is also sufficient for an open neighborhood
of our sequence $S$. Since $\mathcal{S}$ is compact,
this means that we can choose one $\delta$ that suffices
for all sequences, in particular for the sequences
$P_j, P_{j+1}, \ldots$, which completes the proof. \end{proof}

Fix $\delta$ as in the above lemma, and we now fix $l$ such that
$4\tau 2^{-l}<\epsilon$, where $\tau$ is the maximum possible
algebraic degree of the sets $V_n$, the critical values of $P_n$.
Let $\gamma$ be the maximum possible degrees of the algebraic sets
$V_{l,n}$. We can choose $\tau$ and $\gamma$ since the degrees of
the polynomials $P_n$ are bounded from above, which follows from
the compactness of $\mathcal{P}$.

We shall call a holomorphic disc $\Delta$ in a complex line $L$
$\delta$-extendable if the $\delta$-neighborhood of $\Delta$
in $L$ is simply connected.

\begin{lemma}
There exists a constant $c \in \RR$ such that for every $n$ large
enough, every complex line $L$ and every
$\frac{\delta}{4\gamma}$-extendable holomorphic disc
$\Delta\subset L$ that does not intersect a
$\frac{\delta}{2\gamma}$-neighborhood of $L \cap V_{l,n}$, there
exist at least $(1-\epsilon)d(n)^k$ inverse branches of  $P(n)$ on
$\Delta$ for which the preimages $\Delta_i = P(n)^{-1}_i(\Delta)$
satisfy
$$
\mathrm{diam}\left(\Delta_i\right) < cd(n)^{k/2}.
$$
\end{lemma}

\begin{proof}
We can exactly follow the proof of the lemma in [BD] to get that
for every such disc $\Delta$, there exists a constant $c$ such
that there are at least $(1-\epsilon)d^n$ preimages $\Delta_i$ of
diameter less than $cd^{n/2}$. To see that we can choose $c$
independently of $\Delta$, note that we can take the larger disc
$\tilde\Delta$ in that proof as the
$\frac{\delta}{4\gamma}$-neighborhood of $\Delta$ in $L$. It
follows that $\text{Mod}(\tilde\Delta - \Delta)$ is bounded from
below by some strictly positive constant, and this gives a bound
on $c$ which completes the proof of the lemma.
\end{proof}

Note that for every line $L$ that intersects $V(l,n)$ in a finite
number of points and every $x,y$ in the complement of the $\delta
/ \gamma$-neighborhood of $V(l,n)$  in $L$, we can choose a
$\delta/(4\gamma)$-extendable holomorphic disc outside of the
$\delta/(2\gamma)$-neighborhood of $V(l,n)$. Indeed, we can take
the shortest curve in $L$ from $x$ to $y$ that avoids the
$3\delta/(4\gamma)$-neighborhood of $V(l,n)$ and take the
$\delta/(4\gamma)$-neighborhood of the curve as our extendable
disc.

\vskip.25cm

\textbf{Proof of Proposition \ref{preimage}:} Let $X_n$ be the
$\delta$-neighborhood of $V_{l,n}$. We have that $\mu_n(X_n) <
\epsilon$ for any $n \in \NN$. Let $x, y$ be points outside of
$X_n$. We can choose $z$ outside of $X_n$ such that the lines
$L_1$ and $L_2$ through respectively $x, z$ and $y, z$ intersect
$V_{l,n}$ in at most $\gamma$ points. This means that there exist
$\delta/(4\gamma)$-extendable holomorphic discs $\Delta_1 \subset
L_1$ and $\Delta_2 \subset L_2$ such that $x, z \in \Delta_1$ and
$y, z \in \Delta_2$, and such that $\Delta_1$ and $\Delta_2$ avoid
the $\frac{\delta}{2 \gamma}$ neighborhood of $V_{l,n}$. Now it
follows from the lemma that there are at least $(1-\epsilon)d^n$
preimages $x_{j}^{-n}$,  $y_{j}^{-n}$ and $z_j^{-n}$ such that

$$
\text{dist}(x_{j}^{-n}, y_{j}^{-n}) \le \mathrm{dist}(x_{j}^{-n}, z_j^{-n})+ \text{dist}(y_{j}^{-n}, z_j^{-n})
\le 2 \frac{c}{d(n)^{k/2}}.
$$

Hence, for any continuous function $\phi$ of norm $1$ we have that
$$
|\int\phi d\eta_{x,n} - \int\phi d\eta_{y,n}| \le 2\epsilon +
\big| \frac{1}{d^n}\sum_{j} \left( \phi(y_j^{-n})-\phi(y_j^{-n}) \right) \big|
\le 2\epsilon + 2\frac{c}{d(n)^{k/2}}.
$$

For $n$ large enough, this is smaller than $3\epsilon$, which
completes the proof. \hfill $\square$

\vskip.25cm

Now, for some fixed small $\epsilon>0$, let $\epsilon_1,
\epsilon_2, \ldots$ be a monotone decreasing sequence such that
the sum over all $\epsilon_j$ is smaller than $\epsilon$. For
every j, define a set $X_{n,j}$ as in Proposition \ref{preimage} and
$N_j$ in $\NN$ such that $\mu(X_{n,j})<\epsilon_j$ and
$d(\eta_{n,x},\eta_{n,y}) < \epsilon_j$ for any $n$ larger than
$N_j$ and $x,y$ outside of $X_{n,j}$. Now set
$$
U_n := \PP^k - \bigcup_{N_j\le n} X_{n,j}.
$$
We see in particular that $\mu_n(U_n)$ is larger than $1-\epsilon$
for every $n$. Fixing a sequence $x_1, x_2, \ldots$ such that $x_n$ is
an element of $U_n$, we get the following uniform version of Theorem \ref{brolin}.

\begin{lemma} \label{preimages}
For every $\epsilon >0$ there exists an $N$ so that for every $m$ and every $n \ge N$
we have that
$$d(\eta_{x_n, n-m},\mu_m)< \epsilon$$
\end{lemma}

\begin{proof}
We have that
$$
\mu_m = \int \delta_y \, d\mu_m(y),
$$
and therefore we have
$$
\mu_m = \frac{P(m, n+m)^*\mu_{n+m}}{d(m, n + m)^{k}} = \int \eta_{y,n,m} \, d\mu_{n+m}(y).
$$
It follows that
$$
\mu_m - \eta_{x_{n+m},n,m} = \int (\eta_{y,n,m} - \eta_{x_n,n,m}) \, d\mu_n(y).
$$
We can choose a $j$ such that $2\epsilon_j
< \epsilon$, and by our construction of $X_{n,j}$ and $U_n$, it
follows that for $n \ge N_j$ we have $d(\eta_{y,n,m}, \eta_{x_n,n,m})
< \epsilon_j$ for any $y$ outside of $X_{n,j}$, while
$\mu_{n+m}(X_{n+m,j})<\epsilon_j$. Therefore, $d(\mu_m,
\eta_{x_n,j})<2\epsilon_j$, which completes the proof.
\end{proof}

In the autonomous setting it is known that the equilibrium measure is the only
totally invariant measure that doesn't charge the exceptional set \cite{bd}.
We can't expect such a result to hold here.
Consider for instance the map $z \mapsto z^2$ in $\PP^1$.
The equilibrium measures $\mu_n$ are all equal to the normalized Lebesque measure
on the unit circle. However, let $\nu_n$ be the normalized Lebesque measure on the
disc of radius $1/2^n$. Then $\{\nu_n\}$ is totally invariant and doesn't charge the exceptional set $\{0, \infty\}$.

We do have the following related uniqueness result:

\begin{corollary}
Let $\nu$ be a probability measure on $\PP^k$ that doesn't charge
locally pluripolar sets. Then we have that
$$
\frac{P(m, n+m)^*}{d(m, n+m)^k} \nu \rightarrow \mu_m,
$$
weakly.
\end{corollary}

The corollary follows from Proposition \ref{preimage} as in the
proof of lemma \ref{preimages}

\section{Proof of Theorems 1 and 2}

\textbf{Proof of Theorem 1:}
Let $A_0, A_1, \ldots$ be a sequence of measurable subsets
of $\PP^k$ such that $P_n^{-1}(A_n)=A_{n-1}$ for all $n$,
and assume that $\mu_0(A_0)$ is not equal to $0$. We need
to show that $A_0$ has full measure. Define the measures $\nu_n$ by
$$
\nu_n(X)= \frac{\mu_n(X\cap A_n)}{\mu_n(A_n)}.
$$
Clearly, every $\nu_n$ is a probability measure. We see that
\begin{align*}
P_{n*}\nu_{n-1}(X) &= \nu_{n-1}(P_n^{-1}(X))\\
&= \mu_{n-1}(P_n^{-1}(X)\cap A_{n-1}) / \mu_{n-1}(A_{n-1})\\
&= \mu_{n-1}(P_n^{-1}(X \cap A_n)) / \mu_n(A_n)\\
&= P_{n*}\mu_{n-1}(X \cap A_n)/ \mu_n(A_n)\\
&= \mu_n(X\cap A_n)/ \mu_n(A_n) = \nu_n(X).
\end{align*}
Similarly, it follows from the total invariance of the sets $A_n$ and the
measures $\mu_n$ that
$$
\frac{P_ n^*\nu_{n}}{d(n)^k}=\nu_{n-1}.
$$

As we have seen before in the proof of Lemma \ref{preimages}, we have the equation
$$
\mu_0 = \int \eta_{x,n} d\mu_n(x),
$$
and similarly,
$$
\nu_0 = \int \eta_{y,n} d\nu_n(y).
$$
Therefore we see that
$$
\mu_0 - \nu_0 = \int\int (\eta_{x,n} - \eta_{y,n})d\mu_n(x)
\otimes d\nu_n(y).
$$
It now follows from Proposition 1 that for any $\epsilon >0$, we have
$$
\| \mu_0 - \nu_0 \| < 3\epsilon.
$$
Thus $\nu_0 = \mu_0$ and $\mu_0(A_0)$ must equal $1$, which completes the theorem. \hfill $\square$
\vskip.25cm

The argument of the proof of Theorem 2 is similar to that of Theorem 17.1 in [Br].

\vskip.25cm

\textbf{Proof of Theorem 2:}
Let $\phi, \psi$ be test functions of norm at most $1$, and let $\epsilon > 0$.
Construct sets $U_n$ as we did for Lemma \ref{preimages} such that
$\mu_n(U_n) > 1-\epsilon$ for each $n$. It follows from
Lemma \ref{preimages} that we can fix $n$ so large that
$\|\eta_{\zeta,n}-\mu_0\|<\epsilon$ for any $\zeta \in U_n$.

Let $m$ be large enough so that
$$
\int (\phi\circ P(n)) \cdot \psi \, d\mu_0= \int (\phi\circ P(n))\circ \psi
\, d\eta_{x_{m+n},-(m+n)} + \epsilon_1,
$$
where $|\epsilon_1|<\epsilon$. It follows from the definition of
$\eta_{x_{m+n},-(m+n)}$ that the right hand side is equal to
\begin{align*}
\sum_\nu \phi(P(n)(\zeta^\nu_{m+n,-(m+n)}))\psi(\zeta^\nu_{m+n,-(m+n)})
d(m+n)^{-k}+\epsilon_1 \\
= \sum_{\sigma} \phi(\zeta^\sigma_{m+n,-m})d(n, n+m)^{-k}
\sum_{\zeta^\sigma_{m+n,-m} \text{fixed}}
\psi(\zeta^\nu_{m+n,-(m+n)}) d(n)^{-k} +\epsilon_1.
\end{align*}
Counting multiplicity, there are $d(n, n+m)^{k}$ preimages
$\zeta^\sigma_{m+n,-m}$, and since \linebreak
$\mu_n(U_n) > 1-\epsilon$,
we can increase $m$ if necessary so that
 at least $(1-\epsilon)d(n, n+m)^{k}$ of the $\zeta_{m+n,-m}$ are in $U_n$.
It follows that the above right hand side is equal to
$$
\sum \phi(\zeta_{m+n,-m})d(m, n+m)^{-k} \left( \int \psi \, d\mu_0
+\epsilon_3 \right) +\epsilon_1 + \epsilon_2,
$$
where $\epsilon_3$, which depends on $\nu$, and $\epsilon_2$
all have absolute value less than $\epsilon$. We can rewrite this as
$$
\left(\int \psi \, d\mu_0 +\epsilon_3\right) \sum
\phi(\zeta_{m+n,-m})d(n, n+m)^{-k}  +\epsilon_1 + \epsilon_2,
$$
where $\epsilon_3$ no longer depends on $m$. By increasing $m$ if necessary we get
$$
\left(\int \psi \, d\mu_0 +\epsilon_3\right) \left(\int
\phi \, d\mu_n + \epsilon_4 \right) +\epsilon_1 + \epsilon_2,
$$
and so
$$
|\int (\phi \circ P(n)) \cdot \psi \, d\mu_0 - \int \phi \,
d\mu_n \int \psi \, d\mu_0| < 4\epsilon.
$$
This proves that
$$
\int (\phi\circ P(n)) \cdot \psi \, d\mu_0 - \int \phi \, d\mu_n \int \psi \, d\mu_0 \rightarrow 0
$$
for all test functions $\phi$ and $\psi$. The theorem follows
since we can uniformly approximate any continuous function by test
functions.
\hfill $\square$

\begin{remark}
It is not clear if the theorem holds if we allow $\phi$
in the definition of randomly mixing to be in the intersection
of all $\mathcal{L}^2(\mu_n)$, since in general we will not
be able to approximate these functions
by continuous functions that are close in every $\mathcal{L}^2(\mu_n)$
norm at the same time. The theorem does however hold for
$\psi$ in $\mathcal{L}^2(\mu_0)$.
\end{remark}

\section{Random ergodicity in the autonomous setting}

We have already seen that random ergodicity is not equivalent to
ergodicity in the classical case. Indeed, an automorphism can
never have interesting measures that are randomly ergodic, so we
can not expect randomly ergodic to be equivalent to any known
condition from ergodic theory. We shall see in this section
that random ergodicity implies a condition that is stronger
than ergodicity, namely weakly mixing. Recall that a measure
preserving transformation $(P,\mu)$ is \emph{weakly mixing} if for all
$\phi, \psi \in \mathcal{L}^2(\mu)$ we have that
$$
 \frac{1}{n} \sum_{k=1}^{n} \left( \int (\phi\circ P^k) \cdot \psi
\, d\mu - \int \phi \, d\mu \int \psi \, d\mu\right)^2,
$$
converges to $0$ as $n \to \infty$. Weakly mixing implies ergodicity,
see for example \cite{cfs}.

Let $U_P$ be the adjoint
operator working on $\mathcal{L}^2(\mu)$, that is,
$U_P(\phi)=\phi\circ P$.

The following result is well known and can be found in \cite{cfs}

\begin{theorem}
A measure preserving transformation $P$ is weakly mixing if and
only if every eigenfunction of $U_P$ is constant almost
everywhere.
\end{theorem}

We can similarly express random ergodicity in terms of the
operators $U_{P_n}$.

\begin{lemma}
$\{P_n, \mu_n\}$ is randomly ergodic if and only if we have the following
property:

For every sequence $f_0, f_1, \ldots$ with $f_n \in \mathcal{L}^2(\mu_n)$ for
which $f_n\circ P_n = f_{n-1}$ holds for every $n$, we have that
$f_0$ is constant a.e..
\end{lemma}
\begin{proof}
First, assume that the system is randomly ergodic. Fix a totally
invariant sequence of maps $f_0, f_1, \ldots$ as above. For some
$r \in \RR$, define $A_k = \{z | f_k(z) > r\}$. Then $x\in A_{k-1}$
if and only if $f_k(P_k(z)) > r$, and thus if and only if $P_k(z)
\in A_k$, i.e. $P_k^{-1}A_k = A_{k-1}$. This means that $\mu(A_0)
= 0$ or $1$, and this holds for every $r \in \RR$ and thus we see
that the function $f_0$ is constant.

For the converse, let $A_0, \ldots$ be such that $P_k^{-1}(A_k)=
A_{k-1}$ and define $f_k = \mathbf{1}_{A_k}$. It follows that
$f_0$ is constant, and therefore that $A_0$ is has mass $0$ or
$1$.
\end{proof}

\begin{proposition}
Any randomly ergodic measure preserving transformation $(P,\mu)$ is weakly mixing.
\end{proposition}
\begin{proof}
Suppose that $U_P$ has an eigenfunction $f$, say $f\circ P= \sigma
f$. Then define $f_0 = f, f_1= \sigma^{-1}f,$ and so forth. Clearly
this sequence satisfies $f_k \circ P = f_{k-1}$, and therefore
$f_0=f$ is constant a.e..
\end{proof}

\end{document}